\documentclass[final,sort&compress,3p,times]{elsarticle}
\usepackage[]{amsmath}
\usepackage{mathrsfs}
\usepackage{multirow}
\begin{document}
\begin{frontmatter}

\title{The construction of two-dimensional optimal systems for the invariant solutions}

\author[label1]{Xiaorui Hu \corref{cor1}}
\ead{lansexiaoer@163.com}
\author[label2]{Yuqi Li}
\author[label2]{Yong Chen\corref{cor1}}
\ead{ychen@sei.ecnu.edu.cn}

\address[label1]{Department of Applied Mathematics, Zhejiang University of Technology, Hangzhou 310023, China}

\address[label2]{Shanghai Key Laboratory of Trustworthy Computing,
East China Normal University, Shanghai 200062, China}

\cortext[cor1]{Corresponding author.}

\begin{abstract}

To search for inequivalent group invariant solutions, a general and systematic approach is  established to construct two-dimensional optimal systems, which is based on commutator relations, adjoint matrix and the invariants. The details of computing all the invariants for two-dimensional subalgebras is presented and the
optimality of two-dimensional optimal systems is shown clearly
under different values of invariants, with no further proof.
Applying  the algorithm to (1+1)-dimensional heat equation and (2+1)-dimensional Navier-Stokes (NS) equation, their
two-dimensional optimal systems are obtained, respectively. For the heat equation, eleven two-parameter elements in the optimal system are found one by one, which are discovered  more comprehensive. The two-dimensional optimal system of NS equations is used to generate intrinsically different reduced ordinary differential equations and some interesting explicit solutions
are provided.

\end{abstract}


\end{frontmatter}

\section{Introduction}

The study of group-invariant solutions of differential
equations plays an important role in mathematics and physics.
The machinery of Lie group theory provides the systematic
method to search for these special group invariant solutions.
For a $n$-dimensional differential system,
any of its $m$-dimensional ($m<n$) symmetry subgroup can transform
it into a $(n-m)$-dimensional differential system, which is
generally easier to solve than the original system. By solving these
reduced equations, rich  group-invariant solutions are found. For
two  group-invariant solutions, one  may connect them with some group transformation and in this case, one calls them equivalent. Naturally, it is  a significant job to find these inequivalent branches of group-invariant solutions, which leads to the concept of the
optimal systems. For the classification of group-invariant solutions,
it is more convenient to work in the space of the Lie
algebra and this problem reduces to the problem of finding an optimal
system of subalgebras under the adjoint representation .

The adjoint representation of a Lie group on its Lie algebra was known to Lie. The construction of the one-dimensional optimal system of Lie algebra is demonstrated by Ovsiannikov~\cite{Ovsiannikov}, using a global matrix for the adjoint transformation. This is also the technique used by Galas \cite{Galas} and Ibragimov \cite{Ibragimov}. Then Olver~\cite{Olver} uses a slightly different and elegant technique for one-dimensional optimal system, which is based on commutator table and adjoint table, and presents detailed instructions on the KdV equation and the heat equation. For the two-dimensional optimal systems, Ovsiannikov sketches the construction by  showing a simple example. Galas refines Ovsiannikov's method  by  removing equivalent subalgebras for the solvable algebra, and he also discusses the problem of a nonsolvable algebra, which is generally harder. In Ref.~\cite{Coggeshalla}, the details of the construction of the two-dimensional optimal systems are shown for the three-dimensional, one-temperature hydrodynamic equations.
In a fundamental series of papers, Patera et a1. ~\cite{Patera1,Patera2,Patera3,Patera4} have developed a different and powerful method to classify subalgebras
and many  optimal systems of important Lie algebras arising in mathematical physics are obtained.

For the construction of one-dimensional optimal system, Olver has
pointed out that the Killing form of the Lie algebra as an ``invariant" for the adjoint representation is very important since it places restrictions on how far one can expect to simplify the Lie algebra. Recently, Chou et al. ~\cite{Qu1,Qu2,Qu3} develop this idea by introducing more numerical invariants (which are different from common invariants such as the Casmir operator, harmonics and rational invariants) to demonstrate the inequivalence among the elements in the optimal system.
However, to the best of our knowledge,
in spite of the importance of the common invariants for the Lie algebra, there are few literatures to use more invariants except the killing form in the process of constructing optimal systems. For this, in our early paper, we introduce a direct and valid method for providing all the general invariants which are not  numerical invariants and then make the best of them with the adjoint matrix to construct one-dimensional optimal system. On the basis of
all the invariants, the new method can both guarantee the comprehensiveness and the inequivalence of the one-dimensional optimal
system. In this paper, we continue to introduce a systematic method
for two-dimensional optimal system, taking full advantage of the invariants of two-dimensional subalgebras. In almost all of literatures,
one-dimensional optimal system is required for the calculation of
two-dimensional optimal system. Here our new algorithm will start from
Lie algebra directly, without  the prior construction of one-dimensional optimal system. We shall demonstrate the new technique by treating a couple of illustrative examples, the (1+1)-dimensional heat equation and (2+1)-dimensional Navier-Stokes (NS) equation. Due to the two-dimensional optimal system of Lie algebra, the number of independent variables of any differential equations  would be  reduced by two.

The layout of this paper is as follows. In section 2, a systematic algorithm of two-dimensional optimal systems for the
general symmetry algebra is proposed. Since the realization of our new algorithm builds on different invariants and
the adjoint matrix, a valid method for computing all the invariants
of the two-dimensional subalgebras is also given in this section. In section 3, we apply the new algorithm to a six-dimensional Lie algebra of (1+1)-dimensional heat equation, and construct its two-dimensional optimal system step by step. In section 4, the two-dimensional optimal system of
(2+1)-dimensional Navier-Stokes equation is presented and
all the corresponding reduced
equations with some interesting exact group invariant solutions are obtained Finally, a brief conclusion is given in section 5.

\section{A general algorithm for constructing two-dimensional optimal system}\label{algorithm}

Consider the $n$-dimensional symmetry algebra $\mathcal{G}$ of a differential
system, which is generated by the vector fields $\{v_1,v_2, \cdots, v_n\}$. The corresponding symmetry group of $\mathcal{G}$ is denoted as $G$. A family of $r$-dimensional subalgebras $\{\mathfrak{g}_{\alpha}\}_{\alpha \in \mathcal{A}}$ is an $r$-parameter optimal system if
(1) any $r$-dimensional subalgebra is equivalent to some $\mathfrak{g}_{\alpha}$ and (2) $\mathfrak{g}_{\alpha}$ and
$\mathfrak{g}_{\beta}$ are inequivalent for distinct $\alpha$ and $\beta$.

Let
\begin{equation}
\label{Eq1}
w_1\equiv \sum\limits_{i=1} ^{n} a_iv_i, \quad w_2\equiv \sum\limits_{j=1} ^{n} b_jv_j,\quad
w'_1\equiv\sum\limits_{i=1} ^{n} a'_iv_i, \quad w'_2\equiv \sum\limits_{j=1} ^{n} b'_jv_j.
\end{equation}

For two-dimensional optimal system denoted by $\mathfrak{\Theta}_2$, we call its two elements $\mathfrak{g}_{\alpha}=\{w_1,w_2\}$  and $\mathfrak{g}_{\alpha'}=\{w'_1,w'_2\}$ are
equivalent if  one can find
some transformation $g \in G$ and some constants $\{k_1, k_2, k_3, k_4\}$ so that
\begin{equation}
\label{Radw12}
w'_1=k_1Ad_g(w_1)+k_2Ad_g(w_2), \quad w'_2=k_3Ad_g(w_1)+k_4Ad_g(w_2)
\end{equation}
and the inverse is also true. For the invertibility of (\ref{Radw12}),
it requires $k_1k_4-k_2k_3\neq0$.

Here $Ad_g(w_i) (i=1,2)$ is the adjoint representation  action of the group $g$ on the algebra $ w_i$ with $Ad_g(w_i) =g^{-1}w_i g$.

%The main tools used in our
%algorithm are the adjoint transformation matrix and the invariants.
\subsection{Construction of the refined two-dimensional algebra}

For constructing  two-dimensional optimal system $\mathfrak{\Theta}_2$ of $n$-dimensional Lie algebra $\mathcal{G}$, what we should do is to form a list of two-dimensional algebras $\mathcal{G}(w_1,w_2)$ that will later be separated into equivalence classes for different $a_i$ and $b_j$ under the adjoint transformation.
For an arbitrary element $\mathfrak{g}_{\alpha}=\{w_1,w_2\} \in \mathfrak{\Theta}_2$, firstly it requires  $w_1$ and $w_2$ form a
two-dimensional subalgebra, i.e. $[w_1,w_2]=\lambda w_1+\mu w_2$ with
$\lambda$ and $\mu$ being constants. Here $[,]$  represents the commutator
relation. Then Galas refined this selection
by showing that $w_2$ must be an element from the normalizer of $w_1$.
That is to say one can select $w_2$ for
\begin{equation}
\label{RAadw12}
[w_1,w_2]=\lambda w_1.
\end{equation}
Similarly, if $\mathfrak{g}_{\alpha'}=\{w'_1,w'_2\}$ is  equivalent to $\mathfrak{g}_{\alpha}$, it is necessary that
\begin{equation}
\label{RBadw12}
[w'_1,w'_2]=\lambda' w'_1.
\end{equation}

Substituting  (\ref{Eq1}) into (\ref{RAadw12}) and collecting all the coefficients of $v_i$, it will give some relations among $a_1,\cdots a_n, b_1, \cdots b_n$ and $\lambda$, which are called  \emph{the determined equations} by us. What we need do is to find all the representative elements in
$\mathfrak{\Theta}_2$ under the  conditions of \emph{the determined equations}.

In fact, we can split \emph{the determined equations} into two inequivalent classes: $\lambda=0$ and $\lambda \neq 0$. For this inequivalence, we give the following remark.

\textbf{\emph{Remark 1:}} If $\mathfrak{g}_{\alpha}=\{w_1,w_2\}$ and $\mathfrak{g}_{\alpha'}=\{w'_1,w'_2\}$
are equivalent and they satisfy
(\ref{RAadw12}) and (\ref{RBadw12}), we have:

1) when $\lambda=0$, it is natural that $\lambda'=0$;

2) when $\lambda \neq 0$,  there must be $\lambda'\neq 0$ with $k_{2}=0$ in (\ref{Radw12}).

Proof:  Since the equivalence of $\{w_1,w_2\}$ and $\{w'_1,w'_2\}$,
Eqs.(\ref{Radw12}) hold for some $g=e^{v}$ with $v\in \mathcal{G}$. We have
\begin{equation}
\begin{aligned}
\label{five}
&[w'_1,w'_2] =[k_{1}Ad_g(w_1)+k_{2}Ad_g(w_2),k_{3}Ad_g(w_1)+k_{4}Ad_g(w_2)]\\
&=(k_1k_4-k_2k_3)[Ad_g(w_1),Ad_g(w_2)]\\
&=(k_1k_4-k_2k_3)[e^{-v} w_1 e^{v}, e^{-v} w_2 e^{v}]\\
&=(k_1k_4-k_2k_3)e^{-v}[w_1,w_2]e^{v}\\
&=(k_1k_4-k_2k_3)e^{-v}(\lambda w_1)e^{v}\\
&=(k_1k_4-k_2k_3)\lambda Ad_g(w_1)
\end{aligned}
\end{equation}
If there is $\lambda=0$, we obtain $[w'_1,w'_2]=0$, i.e. $\lambda'=0$.
If there is $\lambda\neq0$, in the condition of
\begin{equation}
[w'_1,w'_2]=\lambda'w'_1=\lambda'(k_{1}Ad_g(w_1)+k_{2}Ad_g(w_2)),
\end{equation}
there must be $\lambda' \neq0$ and $k_{2}=0$ for $k_1k_4-k_2k_3\neq0$.

\subsection{Calculation  of  the  adjoint  transformation  matrix} \label{matrix}

For $w_1=\sum\limits_{i=1} ^{n} a_iv_i$, its general adjoint  transformation  matrix $A$  is  the  product  of  the  matrices  of  the  separate  adjoint  actions  $A_1, A_2, \cdots, A_n$, each corresponding to $Ad_{exp(\epsilon v_i)}(w_1), i=1\cdots n$.

For example, applying the adjoint action of $v_1$  to $w_1=
\sum\limits_{i=1} ^{n} a_iv_i$ and with the help of adjoint representation table, one has
\begin{equation}
\label{nn}
\begin{aligned}
&Ad_{exp(\epsilon_1 v_1)}(a_1v_1+a_2v_2+\cdots+a_nv_n)\\
&=a_1Ad_{exp(\epsilon_1 v_1)}v_1+a_2Ad_{exp(\epsilon_1 v_1)}v_2+\cdots+a_nAd_{exp(\epsilon_1 v_1)}v_n\\
&=R_1v_1+R_2v_2+\cdots+R_nv_n,
\end{aligned}
\end{equation}
with $R_i\equiv R_i(a_1,a_2,\cdots, a_n,\epsilon_1), i=1\cdots n$.
To be intuitive,  the formula (\ref{nn}) can be rewritten into the following matrix form:
$$v \doteq (a_1,a_2,\cdots, a_n) \xrightarrow
{Ad_{exp(\epsilon_1 v_1)}} (R_1,R_2,\cdots, R_n)=(a_1,a_2,\cdots, a_n)A_1.
$$
Similarly, the  matrices  $A_2, A_3,\cdots A_n$ of  the  separate  adjoint actions of $v_2, v_3,\cdots, v_n$ can be constructed, respectively. Then the general adjoint  transformation  matrix $A$ is the  product  of  $A_1, \cdots, A_n$  taken  in  any  order
\begin{equation}
A\equiv A(\epsilon_1, \epsilon_2, \cdots, \epsilon_n)=A_1A_2\cdots A_n.
\end{equation}
That is to say, applying the most general adjoint action $Ad_g=Ad_{exp(\epsilon_n v_n)}\cdots Ad_{exp(\epsilon_2 v_2)}Ad_{exp(\epsilon_1 v_1)}
$ to $w_1$,  we have
\begin{equation}
\label{solveeq}
w_1 \doteq (a_1,a_2,\cdots, a_n) \xrightarrow{Ad} Ad_g(w_1)\doteq(a_1,a_2,\cdots, a_n)A.
\end{equation}
Likewise, there is
\begin{equation}
\label{solveeq}
w_2 \doteq (b_1,b_2,\cdots, b_n) \xrightarrow{Ad} Ad_g(w_2)\doteq(b_1,b_2,\cdots, b_n)A.
\end{equation}
When $\mathfrak{g}_{\alpha'}=\{w'_1,w'_2\}$ is equivalent to $\{w_1,w_2\}$,
we can rewrite  (\ref{Radw12}) as
\begin{equation}
\label{RReqs}
\left\{
\begin{aligned}
(a'_1,a'_2,\cdots, a'_n)=k_{1}(a_1,a_2,\cdots, a_n)A+k_{2}
(b_1,b_2,\cdots, b_n)A,\\
(b'_1,b'_2,\cdots, b'_n)=k_{3}(a_1,a_2,\cdots, a_n)A+k_{4}
(b_1,b_2,\cdots, b_n)A,\\
\end{aligned}
 \right.
\end{equation}

\textbf{\emph{Remark 2:}} In fact,  Eqs.(\ref{RReqs}) can be regarded as $2n$ algebraic
equations with respect to $\epsilon_1,\ldots, \epsilon_n$ and $k_1, k_2, k_3, k_4$, which will be taken to judge whether any two given two-dimensional algebras
$\{w_1,w_2\}$ and $\{w'_1,w'_2\}$ are equivalent.
If Eqs.(\ref{RReqs}) have the solution, it means that $\{\sum\limits_{i=1} ^{n} a_iv_i, \sum\limits_{j=1} ^{n} b_jv_j\}$ is equivalent to {$\{\sum\limits_{i=1} ^{n} a'_iv_i, \sum\limits_{j=1} ^{n} b'_jv_j \}$; If Eqs.(\ref{RReqs}) are incompatible, it shows that they are inequivalent.

\subsection{Calculation  of  the equations about the  invariants} \label{invariants}

For the two-dimensional subalgebra, a real function $\phi$ is called an invariant if $\phi(a_{11}Ad_g(w_1)+a_{12}Ad_g(w_2),a_{21}Ad_g(w_1)+a_{22}Ad_g(w_2))=\phi(w_1,w_2)$ for any two-dimensional subalgebra $\{w_1,w_2\} $ and all $g \in G$ with $a_{11},a_{12}, a_{21}, a_{22}$  being arbitrary constants.
For a general
two-dimensional subalgebra $\{\sum\limits_{i=1} ^{n} a_iv_i, \sum\limits_{j=1} ^{n} b_jv_j\}$ of $\mathcal{G}$, the corresponding
invariant is a $2n$-dimensional function of $(a_1,\cdots, a_n, b_1, \cdots, b_n)$.
Now we will propose a valid method to compute  all the invariants of two-dimensional subalgebras and further make the best use of them to construct two-dimensional optimal system.

Let $v=\sum\limits_{k=1} ^{n} c_kv_k$ be a general element of $\mathcal{G}$.
In conjunction with the commutator table, we have
\begin{equation}
\begin{aligned}
Ad_{g}(w_1)&=Ad_{exp(\epsilon v)}(w_1)\\
&=w_1-\epsilon [v,w_1]+ \frac{1}{2!}\epsilon^2[v,[v,w_1]]-\cdots \label{ad1}
\\
&=(a_1v_1+\cdots+a_nv_n)-\epsilon [c_1v_1+\cdots+c_nv_n,a_1v_1+\cdots+a_nv_n]+o(\epsilon) \\
&=(a_1v_1+\cdots+a_nv_n)-\epsilon(\Theta^{a}_1v_1+\ldots+\Theta^{a}_nv_n)
                      +o(\epsilon) \\
&=(a_1-\epsilon\Theta^{a}_1)v_1+(a_2-\epsilon\Theta^{a}_2)v_2+\cdots+
(a_n-\epsilon\Theta^{a}_n)v_n+o(\epsilon),
\end{aligned}
\end{equation}
where $\Theta^{a}_i \equiv \Theta^{a}_i(a_1,\cdots, a_n, c_1,\cdots, c_n)$ can be easily obtained from the commutator table. Similarly, applying  the same adjoint action $v=\sum\limits_{k=1} ^{n} c_kv_k$ to
$w_2$, we get
\begin{equation}
\begin{aligned}
Ad_{g}(w_2)&=Ad_{exp(\epsilon v)}(w_2)\\
&=(b_1-\epsilon\Theta^{b}_1)v_1+(b_2-\epsilon\Theta^{b}_2)v_2+\cdots+
(b_6-\epsilon\Theta^{b}_n)v_n+o(\epsilon),
\end{aligned}
\end{equation}
where $\Theta^{b}_i \equiv \Theta^{b}_i(b_1,\cdots, b_n, c_1,\cdots, c_n)$ is obtained directly by replacing $a_i$ with $b_i$
in $\Theta^{a}_i (i=1 \cdots n)$.

More intuitively, we denote
\begin{equation}
\begin{aligned}
\label{Ad1234}
w_1\doteq (a_1,a_2, \cdots,a_n),\quad w_2 \doteq (b_1,b_2, \cdots,b_n),\\
Ad_{g}(w_1) \doteq (a_1-\epsilon\Theta^{a}_1, a_2-\epsilon\Theta^{a}_2,\cdots, a_n-\epsilon\Theta^{a}_n)+o(\epsilon),\\
Ad_{g}(w_2) \doteq (b_1-\epsilon\Theta^{b}_1, b_2-\epsilon\Theta^{b}_2,\cdots, b_n-\epsilon\Theta^{b}_n)+o(\epsilon).
\end{aligned}
\end{equation}
For the two-dimensional subalgebra $\{w_1,w_2\}$, according to the definition of the invariant we have
\begin{equation}
\begin{aligned}
\label{Rinvariant}
\phi(a_{11}Ad_g(w_1)+a_{12}Ad_g(w_2),a_{21}Ad_g(w_1)+a_{22}Ad_g(w_2))=\phi(w_1,w_2).
\end{aligned}
\end{equation}
Further, to guarantee
$a_{11}Ad_g(w_1)+a_{12}Ad_g(w_2)=w_1$ and $a_{21}Ad_g(w_1)+a_{22}Ad_g(w_2)=w_2$ after the substitution of $\epsilon=0$, we require
\begin{equation}
\label{a1234}
a_{11}\equiv 1+\epsilon a_{11}, \quad  a_{12}\equiv \epsilon a_{12}, \quad
a_{21}\equiv \epsilon a_{21}, \quad  a_{22}\equiv 1+\epsilon a_{22}.
\end{equation}

\textbf{\emph{Remark 3}}: Following ``\emph{Remark 1}'',  we need consider two cases to determine the invariants $\phi$.

(a) When $[w_1, w_2]=0$ i.e $\lambda=0$, substituting (\ref{Ad1234}) and (\ref{a1234}) into Eq.~(\ref{Rinvariant}), then taking the derivative of Eq.~(\ref{Rinvariant}) with respect to $\epsilon$ and
setting $\epsilon=0$,  extracting all the coefficients of  $c_i,a_{11},a_{12},a_{21}, a_{22}$, some linear differential equations of $\phi$ are obtained. By solving these equations, all the invariants $\phi$ can be found.

(b )When $[w_1, w_2]\neq 0$ i.e $\lambda\neq0$, there must be $a_{12}=0$ in Eq.~(\ref{Rinvariant}).  Then one does the same steps just as case (a) to obtain
linear differential equations of $\phi$.

\subsection{the construction of two-dimensional optimal system}

(1) \emph{First step}: give the commutator table and the adjoint representation table of the generators $\{v_i\}^n_{i=1}$ which are firstly
used by Olver~\cite{Olver}. Then following 2.3 and 2.4, compute the corresponding   adjoint  transformation  matrix $A$ and determine the equations about the invariants $\phi$.

(2) \emph{Second step}: in terms of $[w_1, w_2]=\lambda w_1$, present ``\emph{the
determined equations}''  about $a_1, a_2, \cdots, a_n,b_1, b_2, \cdots, b_n$ and $\lambda$. For these equations, two cases including $\lambda=0$ and $\lambda\neq0$
need to be discussed, respectively. The aim is to find all the inequivalent elements $\{\mathfrak{g}_{\alpha}\}$ which form the two-dimensional optimal system.

(3) \emph{Third step}:
for the general element $\{w_1, w_2\}$,  in terms of every
restricted condition shown in  \emph{the
determined equations}, compute their respective invariants
and select the corresponding eligible
representative elements $\{w'_1, w'_2\}$.  For ease of calculations, we  rewrite (\ref{Radw12}) as
\begin{equation}
\left\{
\begin{aligned}
Ad_{g}(w_1)=k'_1w'_1+k'_2w'_2,\\
Ad_{g}(w_2)=k'_3w'_1+k'_4w'_2.
\end{aligned}
 \right.
\quad \quad (k'_1k'_4\neq k'_2k'_3)
\end{equation}
That is to say
\begin{equation}
\label{Reqs}
\left\{
\begin{aligned}
(a_1,a_2,\cdots, a_n)A=k'_{1}(a'_1,a'_2,\cdots, a'_n)+k'_{2}
(b'_1,b'_2,\cdots, b'_n),\\
(b_1,b_2,\cdots, b_n)A=k'_{3}(a'_1,a'_2,\cdots, a'_n)+k'_{4}
(b'_1,b'_2,\cdots, b'_n).\\
\end{aligned}
 \right.
 \quad \quad (k'_1k'_4\neq k'_2k'_3)
\end{equation}
Following ``\emph{Remark 2}'', if Eqs.(\ref{Reqs}) have the solution with respect to $\epsilon_1, \cdots \epsilon_n, k'_1, k'_2, k'_3, k'_4$, it signifies that the selected representative element is right;
if Eqs.(\ref{Reqs}) have no solution, another new representative element  $\{w''_1,w''_2\}$ should be reselected. Repeat the process until all the cases are finished in the ``\emph{the
determined equations}''.

\section{the new approach for  the (1+1)-dimensional heat equation}
The equation for the conduction of heat in a one-dimensional road is written as
\begin{equation}
u_t=u_{xx}.                         \label{}
\end{equation}
The Lie algebra of infinitesimal symmetries for this equation is
spanned by  six vector fields
\begin{equation}
\label{vector}
\begin{aligned}
&v_1=\partial_{x},\quad  v_2=\partial_{t}, \quad
v_3=u\partial_u,\quad
v_4=x\partial_{x}+2t\partial_t,\\
&v_5=2t\partial_{x}-xu\partial_{u}, \quad
v_6=4tx\partial_{x}+4t^2\partial_{t}-(x^2+2t)u\partial_{u},
\end{aligned}
\end{equation}
and the infinitesimal subalgebra
\begin{equation}
v_{h}=h(x,t)\partial_{u},                    \nonumber
\end{equation}
where $h(x,t)$ is an arbitrary solution of the heat equation.
Since the infinite-dimensional subalgebra $\langle{v_h}\rangle$ does not lead to group invariant solutions, it will not be considered in the classification problem.

The  commutator table and actions of the
adjoint representation, which are taken from [6], are in table 1 and table 2,
respectively.

\begin{table}[htbp]
\centering
\caption{\label{table1}the commutator table; the $(i,j)$-th entry is $[v_i,v_j]=v_iv_j-v_jv_i$}
\begin{tabular}{c|cccccc}
\hline
    &  $v_1$  &  $v_2$  &  $v_3$  &  $v_4$  & $v_5$    &  $v_6$ \\
\hline
$v_1$ &  0    &  0    &  0    &  $v_1$  &  $-v_3$  &  2$v_5$ \\
\hline
$v_2$ &  0    &  0    &  0    &  2$v_2$ &  2$v_1$  &  $4v_4-2v_3$ \\
\hline
$v_3$ &  0    &  0    &  0    &  0    &   0    &  0    \\
\hline
$v_4$ &  $-v_1$   &  $-2v_2$    &  0    &  0    &   $v_5$    &  2$v_6$    \\
\hline
$v_5$ &  $v_3$   &   $-2v_1$    &  0    &  $-v_5$    &   0    &  0    \\
\hline
$v_6$ &  $-2v_5$   &  $2v_3-4v_4$    &  0    &  $-2v_6$    &   0    &  0    \\
\hline
\end{tabular}
\end{table}

\begin{table}[htbp]
\centering
\caption{\label{table2}the adjoint representation
table;the $(i,j)$-th entry gives $Ad_{exp(\epsilon v_i)}(v_j)$}
\begin{tabular}{c|cccccc}
\hline
 Ad   &  $v_1$  &  $v_2$  &  $v_3$  &  $v_4$  & $v_5$    &  $v_6$ \\
\hline
$v_1$ &  $v_1$    &  $v_2$    &  $v_3$    &  $v_4-\epsilon v_1 $  &  $v_5+\epsilon v_3$  &  $v_6-2\epsilon v_5-\epsilon^2 v_3$ \\
\hline
$v_2$ &  $v_1$    &  $v_2$     &  $v_3$     & $v_4-2\epsilon v_2$ &  $v_5-2\epsilon v_1$  &  $v_6-4\epsilon v_4+2\epsilon v_3+4\epsilon^2v_2$ \\
\hline
$v_3$ &  $v_1$    &  $v_2$     &  $v_3$  &   $v_4$    &   $v_5$    &  $v_6$   \\
\hline
$v_4$ &  $e^{\epsilon}v_1$   &  $e^{2\epsilon}v_2$    &  $v_3$    &  $v_4$    &   $e^{-\epsilon}v_5$    &  $e^{-2\epsilon}v_6$    \\
\hline
$v_5$ &  $v_1-\epsilon v_3$   &   $v_2+2\epsilon v_1-\epsilon^2 v_3$    & $v_3$ & $v_4+\epsilon v_5$    &  $v_5$    &  $v_6$    \\
\hline
$v_6$ &  $v_1+2\epsilon v_5$   &  $v_2-2\epsilon v_3+4\epsilon
v_4+4\epsilon^2 v_6$  & $v_3$ &  $v_4+2\epsilon v_6$    &  $v_5$    &  $v_6$  \\
\hline
\end{tabular}
\end{table}

\subsection{Calculation  of  the  adjoint  transformation  matrix} \label{matrix}

For $w_1=\sum\limits_{i=1} ^{6} a_iv_i$, its general adjoint  transformation  matrix $A$  is  the  product  of  the  matrices  of  the  separate  adjoint  actions  $A_1, A_2, \cdots, A_6$, each corresponding to $Ad_{exp(\epsilon v_i)}(w_1), i=1\cdots 6$.

For example, under the adjoint action of $v_1$  and with the help of Table 2,
$w_1$ can be transformed into
\begin{equation}
\label{a1}
\begin{aligned}
&Ad_{exp(\epsilon_1 v_1)}(a_1v_1+a_2v_2+a_3v_3+a_4v_4+a_5v_5+a_6v_6)\\
&=(a_1-a_4\epsilon_1)v_1+a_2v_2+(a_3+a_5\epsilon_1-a_6\epsilon_1^2)v_3+a_4v_4+(a_5-2\epsilon_1 a_6)v_5+a_6v_6.
\end{aligned}
\end{equation}
One can rewrite above formula (\ref{a1}) into the following matrix form:
$$w_1 \doteq (a_1,a_2,\cdots, a_6) \xrightarrow{Ad_{exp(\epsilon_1 v_1)}} (a_1,a_2,\cdots, a_6)A_1
$$
where
\begin{equation}
A_1=\left(
\begin{array}{cccccc}
 1 & 0 & 0 & 0 & 0 & 0 \\
 0 & 1 & 0 & 0 & 0 & 0\\
 0 & 0 & 1 & 0 & 0 & 0 \\
 -\epsilon_1 & 0 & 0 & 1 & 0 & 0 \\
 0 & 0 & \epsilon_1 & 0 & 1 & 0 \\
 0 & 0 & \epsilon_1^2 & 0 & -2\epsilon_1 & 1
 \end{array}
\right).
\end{equation}
Similarly,  the  rest matrices  of  the  separate  adjoint actions of
$v_2, \cdots, v_6$  are  found  to  be:

\begin{equation}
A_2=\left(
\begin{array}{cccccc}
 1 & 0 & 0 & 0 & 0 & 0 \\
 0 & 1 & 0 & 0 & 0 & 0\\
 0 & 0 & 1 & 0 & 0 & 0 \\
 0 &-2\epsilon_2 & 0 & 1 & 0 & 0 \\
 -2\epsilon_2 & 0 & 0 & 0 & 1 & 0 \\
 0 & 4\epsilon_2^2 & 2\epsilon_2 & -4\epsilon_2 & 0 & 1
 \end{array}
\right),\quad
A_4=\left(
\begin{array}{cccccc}
 e^{\epsilon_4} & 0 & 0 & 0 & 0 & 0 \\
 0 & e^{2\epsilon_4} & 0 & 0 & 0 & 0\\
 0 & 0 & 1 & 0 & 0 & 0 \\
 0 & 0 & 0 & 1 & 0 & 0 \\
 0 & 0 & 0 & 0 & e^{-\epsilon_4} & 0 \\
 0 & 0 & 0 & 0 & 0 & e^{-2\epsilon_4}
 \end{array}
\right).
\end{equation}

\begin{equation}
A_5=\left(
\begin{array}{cccccc}
 1 & 0 & -\epsilon_5 & 0 & 0 & 0 \\
 2\epsilon_5 & 1 & -\epsilon_5^2 & 0 & 0 & 0\\
 0 & 0 & 1 & 0 & 0 & 0 \\
 0 &0 & 0 & 1 & \epsilon_5 & 0 \\
 0 & 0 & 0 & 0 & 1 & 0 \\
 0 &0 & 0 & 0 & 0 & 1
 \end{array}
\right),\quad \quad  \quad \quad
A_6=\left(
\begin{array}{cccccc}
1  & 0 & 0 & 0 & 2\epsilon_6 & 0 \\
 0 & 1 & -2\epsilon_6 & 4\epsilon_6 & 0 & 4\epsilon_6^2\\
 0 & 0 & 1 & 0 & 0 & 0 \\
 0 & 0 & 0 & 1 & 0 & 2\epsilon_6 \\
 0 & 0 & 0 & 0 & 1 & 0 \\
 0 & 0 & 0 & 0 & 0 & 1
 \end{array}
\right),
\end{equation}
with $A_3=E$ being the identity matrix. Then the general adjoint  transformation  matrix $A$ is the  product  of  $A_1, \cdots, A_6$  which can be taken  in  any  order:
\begin{equation}
\label{AA}
A\equiv(a_{ij})_{6\times6}=A_1A_2A_3A_4A_5A_6,
\end{equation}
with
\begin{equation}
\begin{aligned}
&a_{11}=e^{\epsilon_4}, \quad a_{12}=0, \quad a_{13}=-\epsilon_5e^{\epsilon_4}, \quad a_{14}=0, \quad a_{15}=2\epsilon_6 e^{\epsilon_4}, \quad a_{16}=0, \quad
a_{21}= 2\epsilon_5 e^{2\epsilon_{4}},\quad a_{22}= e^{2\epsilon_4},\\
&a_{23}=-(\epsilon_5^2+2\epsilon_6)e^{2\epsilon_4}, \quad a_{24}=4\epsilon_6e^{2\epsilon_4}, \quad a_{25}=4\epsilon_5\epsilon_6 e^{2\epsilon_4},\quad a_{26}=4\epsilon_6^2 e^{2\epsilon_4},\quad a_{31}=a_{32}=a_{34}=a_{35}=a_{36}=0, \\
&a_{33}=1,
\quad
a_{41}=-e^{\epsilon_4}
(\epsilon_1+4\epsilon_2\epsilon_5e^{\epsilon_4}), \quad a_{42}=-2\epsilon_2e^{2\epsilon_4},\quad a_{43}=4\epsilon_2\epsilon_6e^{2\epsilon_4}+\epsilon_5e^{\epsilon_4}(\epsilon_1
+2\epsilon_2\epsilon_5e^{\epsilon_4}), \\ &a_{44}=1-8\epsilon_2\epsilon_6e^{2\epsilon_4}, \quad a_{45}=\epsilon_5-2\epsilon_6e^{\epsilon_4}(\epsilon_1+4\epsilon_2\epsilon_5e^{\epsilon_4}),\quad
a_{46}=2\epsilon_6
(1-4\epsilon_2\epsilon_6e^{2\epsilon_4}),\quad a_{51}=-2\epsilon_2e^{\epsilon_4}, \quad a_{52}=0,\\
&a_{53}=\epsilon_1
+2\epsilon_2\epsilon_5e^{\epsilon_4}, \quad a_{54}=0, \quad a_{55}=e^{-\epsilon_4}(1-4\epsilon_2\epsilon_6e^{2\epsilon_4}),\quad
a_{56}=0,\quad a_{61}=4\epsilon_2e^{\epsilon_4}(\epsilon_1+2\epsilon_2\epsilon_5e^{\epsilon_4}), \\
& a_{62}=4\epsilon^2_2e^{2\epsilon_4},\quad a_{63}=-
(\epsilon_1+2\epsilon_2\epsilon_5e^{\epsilon_4})^2+2\epsilon_2(1-4\epsilon_2\epsilon_6
e^{2\epsilon_4}), \quad a_{64}=-4\epsilon_2(1-4\epsilon_2\epsilon_6e^{2\epsilon_4}),\\ &a_{65}=8\epsilon_2\epsilon_6e^{\epsilon_4}(\epsilon_1
+2\epsilon_2\epsilon_5e^{\epsilon_4})-4\epsilon_2\epsilon_5-2\epsilon_1e^{-\epsilon_4},\quad
a_{66}=e^{-2\epsilon_4}
(1-4\epsilon_2\epsilon_6e^{2\epsilon_4})^2.
\end{aligned}
\end{equation}
\subsection{Calculation  of  the  invariants}
For the heat equation, take
\begin{equation}
\label{w12}
w_1=\sum\limits_{i=1} ^{6} a_iv_i, \quad w_2=\sum\limits_{j=1} ^{6} b_jv_j.
\end{equation}

For a general
two-dimensional subalgebra $\{w_1, w_2\}$ of the heat equation, the corresponding
invariant $\phi$ is a twelve-dimensional function of $(a_1,\cdots, a_6, b_1, \cdots, b_6)$. Let $v=\sum\limits_{k=1} ^{6} c_kv_k$ be a general element of $\mathcal{G}$,
then in conjunction with the commutator table 1, we have
\begin{equation}
\begin{aligned}
Ad_{g}(w_1)&=Ad_{exp(\epsilon v)}(w_1)\\
&=w_1-\epsilon [v,w_1]+ \frac{1}{2!}\epsilon^2[v,[v,w_1]]-\cdots \label{ad1}
\\
&=(a_1v_1+\cdots+a_6v_6)-\epsilon [c_1v_1+\cdots+c_6v_6,a_1v_1+\cdots+a_6v_6]+o(\epsilon) \\
&=(a_1v_1+\cdots+a_6v_6)-\epsilon(\Theta^{a}_1v_1+\ldots+\Theta^{a}_6v_6)
                      +o(\epsilon) \\
&=(a_1-\epsilon\Theta^{a}_1)v_1+(a_2-\epsilon\Theta^{a}_2)v_2+\cdots+
(a_6-\epsilon\Theta^{a}_6)v_6+o(\epsilon)
\end{aligned}
\end{equation}
with
\begin{equation}
\label{thetaa}
\begin{aligned}
&\Theta^{a}_1=-c_4a_1-2c_5a_2+c_1a_4+2c_2a_5,\quad
\Theta^{a}_2=-2c_4a_2+2c_2a_4,\quad
\Theta^{a}_3=c_5a_1+2c_6a_2-c_1a_5-2c_2a_6\\
&\Theta^{a}_4=-4c_6a_2+4c_2a_6,\quad
\Theta^{a}_5=-2c_6a_1-c_5a_4+c_4a_5+2c_1a_6,\quad
\Theta^{a}_6=-2c_6a_4+2c_4a_6.
\end{aligned}
\end{equation}
Similarly, applying  the same adjoint action $v=\sum\limits_{k=1} ^{6} c_kv_k$ to
$w_2$, we get
\begin{equation}
\begin{aligned}
Ad_{g}(w_2)&=Ad_{exp(\epsilon v)}(w_2)\\
&=(b_1-\epsilon\Theta^{b}_1)v_1+(b_2-\epsilon\Theta^{b}_2)v_2+\cdots+
(b_6-\epsilon\Theta^{b}_6)v_6+o(\epsilon),
\end{aligned}
\end{equation}
with
\begin{equation}
\label{thetab}
\begin{aligned}
&\Theta^{b}_1=-c_4b_1-2c_5b_2+c_1b_4+2c_2b_5,\quad
\Theta^{b}_2=-2c_4b_2+2c_2b_4,\quad
\Theta^{b}_3=c_5b_1+2c_6b_2-c_1b_5-2c_2b_6\\
&\Theta^{b}_4=-4c_6b_2+4c_2b_6,\quad
\Theta^{b}_5=-2c_6b_1-c_5b_4+c_4b_5+2c_1b_6,\quad
\Theta^{b}_6=-2c_6b_4+2c_4b_6.
\end{aligned}
\end{equation}
For the two-dimensional subalgebra $\{w_1,w_2\}$, according to the definition of the invariant, we have
\begin{equation}
\begin{aligned}\label{invariant}
\phi((1+\epsilon a_{11})Ad_g(w_1)+\epsilon a_{12}Ad_g(w_2),\epsilon a_{21}Ad_g(w_1)+(1+\epsilon a_{22})Ad_g(w_2))=\phi(w_1,w_2).
\end{aligned}
\end{equation}

Following ``\emph{Remark 1}'' and ``\emph{Remark 3}'', Eq.(\ref{invariant}) can be separated into
two cases.

(a) When $\lambda=0$,  all the $c_i (i=1\cdots 6), a_{11},a_{12},a_{21}, a_{22}$ in Eq.(\ref{invariant}) are arbitrary. Now
taking the derivative of Eq.~(\ref{invariant}) with respect to $\epsilon$ and then
setting $\epsilon=0$, extracting the coefficients of all $c_i,a_{11},a_{12},a_{21}, a_{22}$, one can directly  obtain nine differential equations about $\phi\equiv \phi(a_1,\cdots, a_6, b_1, \cdots, b_6)$ :
\begin{equation}
\label{8Eqs}
\begin{aligned}
a_1\phi_{a_1}+a_2\phi_{a_2}+a_3\phi_{a_3}+a_4\phi_{a_4}+a_5\phi_{a_5}+a_6\phi_{a_6}=0,\\
a_1\phi_{b_1}+a_2\phi_{b_2}+a_3\phi_{b_3}+a_4\phi_{b_4}+a_5\phi_{b_5}+a_6\phi_{b_6}=0,\\
b_1\phi_{b_1}+b_2\phi_{b_2}+b_3\phi_{b_3}+b_4\phi_{b_4}+b_5\phi_{b_5}+b_6\phi_{b_6}=0,\\
2a_2\phi_{a_1}+2b_2\phi_{b_1}-a_1\phi_{a_3}-b_1\phi_{b_3}+a_4\phi_{a_5}+b_4\phi_{b_5}=0,\\
-a_4\phi_{a_1}-b_4\phi_{b_1}+a_5\phi_{a_3}+b_5\phi_{b_3}-2a_6\phi_{a_5}-2b_6\phi_{b_5}=0,\\
a_1\phi_{a_1}+b_1\phi_{b_1}+2a_2\phi_{a_2}+2b_2\phi_{b_2}-a_5\phi_{a_5}-b_5\phi_{b_5}
-2a_6\phi_{a_6}-2b_6\phi_{b_6}=0,\\
-a_5\phi_{a_1}-b_5\phi_{b_1}-a_4\phi_{a_2}-b_4\phi_{b_2}+a_6\phi_{a_3}+b_6\phi_{b_3}
-2a_6\phi_{a_4}-2b_6\phi_{b_4}=0,\\
-a_2\phi_{a_3}-b_2\phi_{b_3}+2a_2\phi_{a_4}+2b_2\phi_{b_4}+a_1\phi_{a_5}+b_1\phi_{b_5}
+a_4\phi_{a_6}+b_4\phi_{b_6}=0,\\
\end{aligned}
\end{equation}
and
\begin{equation}
b_1\phi_{a_1}+b_2\phi_{a_2}+b_3\phi_{a_3}+b_4\phi_{a_4}+b_5\phi_{a_5}+b_6\phi_{a_6}=0.
\label{1Eq}
\end{equation}
Here the subscripts  indicate  partial derivatives.

(b) When $\lambda \neq 0$, we have  $a_{12}=0$ and all the $c_i (i=1\cdots 6), a_{11},a_{21}, a_{22}$ in Eq.(\ref{invariant}) are arbitrary. Plugging  $a_{12}=0$ into Eq.(\ref{invariant}) and making the same process to
case(a), we get eight equations about $\phi$ which are just Eqs.(\ref{8Eqs}).

\subsection{the two-dimensional optimal system for the heat equation}

Due to the refined Lie algebra $[w_1,w_2]=\lambda w_1$,
we find its determined equations
are
\begin{equation}
\label{fengbi}
\begin{aligned}
& a_1b_4+2a_2b_5-a_4b_1-2a_5b_2=\lambda a_1, \quad
2a_2b_4-2a_4b_2=\lambda a_2, \quad
-a_1b_5-2a_2b_6+a_5b_1+2a_6b_2=\lambda a_3,
\\
& 2a_1b_6+a_4b_5-a_5b_4-2a_6b_1=\lambda a_5,  \quad
 4a_2b_6-4a_6b_2=\lambda a_4, \quad
2a_4b_6-2a_6b_4=\lambda a_6.
\end{aligned}
\end{equation}
Then Eqs.(\ref{fengbi}) are separated into two inequivalent
classes: $\lambda=0$  and $\lambda\neq0$.

\subsubsection{the case of $\lambda=0$ in the determined equations (\ref{fengbi})}

Substituting $\lambda=0$ into Eqs.(\ref{fengbi}), two cases need to be considered.

\textbf{\emph{Case 1:}} Not all $a_2, a_4, a_6, b_2, b_4$ and $b_6$ are zeroes. Without loss of generality, we adopt $a_6\neq0$. In fact, when only one of $a_2, a_4, a_6, b_2, b_4$ and $b_6$ is not zero, one can choose  appropriate adjoint transformation  to transform it into the case $\tilde{a}_6\neq0$. By solving Eqs.(\ref{fengbi}) with $\lambda=0$ and $a_6\neq0$, we have three kinds of solutions:

(i) $a_3, a_4, a_5, a_6, b_3$ and $b_6$ are independent with

\begin{equation}
\label{eqi}
a_1=\frac{1}{2}\frac{a_4a_5}{a_6},\quad a_2=\frac{1}{4}\frac{a^2_4}{a_6},\quad
b_1=\frac{1}{2}\frac{a_4a_5b_6}{a^2_6}, \quad b_2=\frac{1}{4}\frac{a^2_4b_6}{a^2_6},\quad
b_4=\frac{a_4b_6}{a_6}, \quad
b_5=\frac{a_5b_6}{a_6}.
\end{equation}

(ii) $a_3, a_4, a_5, a_6, b_3, b_5 $ and $b_6$ are arbitrary but with
\begin{equation}
\label{eqii}
a_1=\frac{1}{2}\frac{a_4a_5}{a_6},\quad a_2=\frac{1}{4}\frac{a^2_4}{a_6},\quad
b_1=\frac{1}{2}\frac{a_4b_5}{a_6}, \quad b_2=\frac{1}{4}\frac{a^2_4b_6}{a^2_6},\quad
b_4=\frac{a_4b_6}{a_6}, \quad
b_5\neq\frac{a_5b_6}{a_6}.
\end{equation}

(iii) $a_1, a_2, a_3, a_4, a_5, a_6, b_3$ and $b_6$ are arbitrary but with
\begin{equation}
\label{eqiii}
a_1\neq \frac{1}{2}\frac{a_4a_5}{a_6} \quad {\rm{or}} \quad a_2\neq \frac{1}{4}\frac{a^2_4}{a_6},\quad
b_1=\frac{a_1b_6}{a_6}, \quad b_2=\frac{a_2b_6}{a_6},\quad
b_4=\frac{a_4b_6}{a_6}, \quad
b_5=\frac{a_5b_6}{a_6}.
\end{equation}
%For the above three solutions, there are  two remarks:\\
%\textbf{\emph{Remark 3:}} Every solution about itself is closed. That is to say,
%if a two-dimensional subalgebra $\{w_1,w_2\}$ of (\ref{w12}) satisfy Eq.(\ref{eqi}) (Eq.(\ref{eqii}) or Eq.(\ref{eqiii})), so do all its equivalent two-dimensional subalgebra $\{a_{11}Ad_g(w_1)+a_{12}Ad_g(w_2),a_{21}Ad_g(w_1)+a_{22}Ad_g(w_2)\}$.\\
%\textbf{\emph{Remark 4:}}

Substituting the above three conditions into Eqs.(\ref{8Eqs}) and Eq.(\ref{1Eq}), we find that $\phi\equiv constant $, i.e. there is
no invariant. Then for each case, select the corresponding representative element
$\{w'_1,w'_2\}$ and verify whether Eqs.(\ref{Reqs}) have the solution.

For case (i), select a representative element  $\{w'_1,w'_2\}=\{v_6,v_3\}$. Substituting (\ref{eqi}) and $w'_1=v_6, w'_2=v_3$ into Eqs.(\ref{Reqs}), we obtain the solution
\begin{equation}
k'_1=a_6e^{-2\epsilon_4}, \quad k'_2=\frac{4a_3a_6+2a_4a_6+a^2_5}{4a_6},\quad
k'_3=b_6e^{-2\epsilon_4}, \quad k'_4=\frac{4b_3a^2_6+b_6a^2_5+2a_4a_6b_6}{4a^2_6},\quad
\epsilon_1=\frac{a_5}{2a_6}, \quad \epsilon_2=\frac{a_4}{4a_6}.
\end{equation}
Hence case (i) is equivalent to $\{v_6,v_3\}$.

For case (ii), there exist three circumstances in terms of the following expression
\begin{equation}
\label{inv1}
\Lambda_1\equiv 2a_6[a_4(a_5b_6-b_5a_6)^2-2a_6(a_3b_6-b_3a_6)^2-2a_6(a_5b_6-b_5a_6)(a_3b_5-b_3a_5)].
\end{equation}\\
(iia) When $\Lambda_1>0$, case (ii) is equivalent to $\{v_3+v_6,v_5\}$. After the substitution of (\ref{eqii}) with $w'_1=v_3+v_6, w'_2=v_5$,  Eqs.(\ref{Reqs}) hold for
\begin{equation}
\begin{aligned}
&k'_1=\frac{\Lambda_1}{4a_6(a_5b_6-b_5a_6)^2},\quad
k'_2=\frac{a_5(a_5b_6-b_5a_6)+2a_6(a_3b_6-b_3a_6)}{2|a_6|(a_5b_6-b_5a_6)^2}
\sqrt{\Lambda_1},\quad
k'_3=\frac{b_6}{a_6}k_1,\quad \epsilon_2=\frac{a_4}{4a_6},\\
&k'_4=\frac{b_5(a_5b_6-b_5a_6)+2b_6(a_3b_6-b_3a_6)}{2|a_6|(a_5b_6-b_5a_6)^2}
\sqrt{\Lambda_1},\quad
\epsilon_1=-\frac{a_3b_6-b_3a_6}{a_5b_6-b_5a_6},\quad
\epsilon_4=\ln{\frac{2|a_6(a_5b_6-b_5a_6)|}{\sqrt{\Lambda_1}}}.
\end{aligned}
\end{equation}\\
(iib) When $\Lambda_1<0$, case (ii) is equivalent to $\{-v_3+v_6,v_5\}$. The solution for Eqs.(\ref{Reqs}) is
\begin{equation}
\begin{aligned}
&k'_1=-\frac{\Lambda_1}{4a_6(a_5b_6-b_5a_6)^2},\quad
k'_2=\frac{a_5(a_5b_6-b_5a_6)+2a_6(a_3b_6-b_3a_6)}{2|a_6|(a_5b_6-b_5a_6)^2}
\sqrt{-\Lambda_1},\quad
k'_3=\frac{b_6}{a_6}k_1,\quad \epsilon_2=\frac{a_4}{4a_6},\\
&k'_4=\frac{b_5(a_5b_6-b_5a_6)+2b_6(a_3b_6-b_3a_6)}{2|a_6|(a_5b_6-b_5a_6)^2}
\sqrt{-\Lambda_1},\quad
\epsilon_1=-\frac{a_3b_6-b_3a_6}{a_5b_6-b_5a_6},\quad
\epsilon_4=\ln{\frac{2|a_6(a_5b_6-b_5a_6)|}{\sqrt{-\Lambda_1}}}.
\end{aligned}
\end{equation}\\
(iic) When $\Lambda_1=0$, case (ii) is equivalent to $\{v_6,v_5\}$. By solving Eqs.(\ref{Reqs}), we obtain
\begin{equation}
\begin{aligned}
&k'_1=a_6e^{-2\epsilon_4},\quad
k'_2=\frac{a_5(a_5b_6-b_5a_6)+2a_6(a_3b_6-b_3a_6)}{(a_5b_6-b_5a_6)e^{\epsilon_4}},\quad k'_4=\frac{a_5(a_5b_6-b_5a_6)+2a_6(a_3b_6-b_3a_6)}{(a_5b_6-b_5a_6)e^{\epsilon_4}},\\
&k'_3=b_6e^{-2\epsilon_4},\quad
\epsilon_1=-\frac{a_3b_6-b_3a_6}{a_5b_6-b_5a_6},\quad
\epsilon_2=\frac{(a_5b_6-b_5a_6)(a_3b_5-a_5b_3)+(a_3b_6-b_3a_6)^2}{2(a_5b_6-b_5a_6)^2}.
\end{aligned}
\end{equation}

For case (iii), it can be divided into the following several  types.\\
(iiia)  $4a_2a_6-a^2_4>0.$ Select a representative element $\{v_2+v_6,v_3\}$. After substituting (\ref{eqiii}) into Eqs.(\ref{Reqs}), we have
\begin{equation}
\begin{aligned}
&k'_1=(a_2-2a_4\epsilon_2+4a_6\epsilon^2_2)e^{2\epsilon_4},\quad
k'_2=a_3+\frac{a_4}{2}+\frac{a^2_1a_6-a_1a_4a_5+a_2a^2_5}{4a_2a_6-a^2_4},\quad
k'_3=\frac{b_6(a_2-2a_4\epsilon_2+4a_6\epsilon^2_2)}{a_6}e^{2\epsilon_4},\\
&k'_4=b_3+\frac{b_6}{a_6}(\frac{a_4}{2}+\frac{a^2_1a_6-a_1a_4a_5+a_2a^2_5}{4a_2a_6-a^2_4}),\quad
\quad \epsilon_1=\frac{2\epsilon_2(2a_1a_6-a_4a_5)+2a_2a_5-a_1a_4}{4a_2a_6-a^2_4},
\\
&\epsilon_5=\frac{a_4a_5-2a_1a_6}{(4a_2a_6-a^2_4)e^{\epsilon_4}},\quad
\epsilon_6=\frac{4a_6\epsilon_2-a_4}{4e^{2\epsilon_4}(a_2-2a_4\epsilon_2
+4a_6\epsilon^2_2)}, \quad e^{2\epsilon_4}=\frac{\sqrt{4a_2a_6-a^2_4}}{2|a_2-2a_4\epsilon_2+4a_6\epsilon^2_2|}.
\end{aligned}
\end{equation}
(iiib)  $4a_2a_6-a^2_4<0.$ Choose a representative element $\{-v_2+v_6,v_3\}$. Now
Eqs.(\ref{Reqs}) have the solution
\begin{equation}
\begin{aligned}
&k'_1=-(a_2-2a_4\epsilon_2+4a_6\epsilon^2_2)e^{2\epsilon_4},\quad
k'_2=a_3+\frac{a_4}{2}+\frac{a^2_1a_6-a_1a_4a_5+a_2a^2_5}{4a_2a_6-a^2_4},\quad
k'_3=-\frac{b_6(a_2-2a_4\epsilon_2+4a_6\epsilon^2_2)}{a_6}e^{2\epsilon_4},\\
&k'_4=b_3+\frac{b_6}{a_6}(\frac{a_4}{2}+\frac{a^2_1a_6-a_1a_4a_5+a_2a^2_5}{4a_2a_6-a^2_4}),\quad
\quad \epsilon_1=\frac{2\epsilon_2(2a_1a_6-a_4a_5)+2a_2a_5-a_1a_4}{4a_2a_6-a^2_4},
\\
&\epsilon_5=\frac{a_4a_5-2a_1a_6}{(4a_2a_6-a^2_4)e^{\epsilon_4}},\quad
\epsilon_6=\frac{4a_6\epsilon_2-a_4}{4e^{2\epsilon_4}(a_2-2a_4\epsilon_2
+4a_6\epsilon^2_2)}, \quad e^{2\epsilon_4}=\frac{\sqrt{-(4a_2a_6-a^2_4)}}{2|a_2-2a_4\epsilon_2+4a_6\epsilon^2_2|}.
\end{aligned}
\end{equation}
(iiic)  $4a_2a_6-a^2_4=0.$ In this case, two conditions should be considered.

When $2a_1a_6-a_4a_5>0$, adopt the representative element $\{v_1+v_6,v_3\}$. Then the solution for Eqs.(\ref{Reqs}) is
\begin{equation}
\begin{aligned}
&k'_1=\frac{a_6}{Z^2},\quad k'_2=-\frac{(\epsilon^2_6+\epsilon_5)(2a_1a_6-a_4a_5)}{2a_6}Z+a_3+\frac{a_4}{2}
+\frac{a^2_5}{4a_6},\quad
k'_3=\frac{b_6(2a_1a_6-a_4a_5)}{2a^2_6}Z,\quad \epsilon_2=\frac{a_4}{4a_6},\\
&k'_4=\frac{b_6}{4a^2_6}[-2(\epsilon_6^2+\epsilon_5)(2a_1a_6-a_4a_5)Z+(a^2_5+2a_4a_6)]
+b_3,\quad
\epsilon_1=\frac{\epsilon_6(2a_1a_6-a_4a_5)}{2a^2_6}Z^2+\frac{a_5}{2a_6},\quad
\epsilon_4=\ln(Z).
\end{aligned}
\end{equation}
with $Z=(\frac{2a^2_6}{2a_1a_6-a_4a_5})^{1/3}$. \\

When $2a_1a_6-a_4a_5<0$, adopt the representative element $\{-v_1+v_6,v_3\}$.
By solving Eqs.(\ref{Reqs}), we find
\begin{equation}
\begin{aligned}
&k'_1=\frac{a_6}{Z'^2},\quad k'_2=-\frac{(\epsilon_5-\epsilon^2_6)(2a_1a_6-a_4a_5)}{2a_6}Z'+a_3+\frac{a_4}{2}
+\frac{a^2_5}{4a_6},\quad
k'_3=-\frac{b_6(2a_1a_6-a_4a_5)}{2a^2_6}Z',\quad \epsilon_2=\frac{a_4}{4a_6},\\
&k'_4=\frac{b_6}{4a^2_6}[-2(\epsilon_5-\epsilon^2_6)(2a_1a_6-a_4a_5)Z'+(a^2_5+2a_4a_6)]
+b_3,\quad
\epsilon_1=\frac{\epsilon_6(2a_1a_6-a_4a_5)}{2a^2_6}Z'^2+\frac{a_5}{2a_6},\quad
\epsilon_4=\ln(Z').
\end{aligned}
\end{equation}
with $Z'=(-\frac{2a^2_6}{2a_1a_6-a_4a_5})^{1/3}$. \\
\textbf{\emph{Case 2:}} $a_2=a_4=a_6=b_2=b_4=b_6=0$. Now the determined equations
(\ref{fengbi}) become
\begin{equation}
\label{case21}
-a_1b_5+a_5b_1=0.
\end{equation}
Here we just need consider not all $a_1, a_5, b_1$  and $b_5$ are zeroes. Without loss of generality, let $a_5\neq0$. Similarly, if one of $a_1, a_5, b_1, b_5$ is not zero, one can choose  appropriate adjoint transformation  to transform it into the case $a_5\neq0$. By solving Eq.(\ref{case21}), we obtain
\begin{equation}
\label{case211}
b_1=\frac{a_1b_5}{a_5}.
\end{equation}
Now adopt a representative element $\{v_5,v_3\}$. Then one can verify that all the
$\{a_1v_1+a_3v_3+a_5v_5, b_1v_1+b_3v_3+b_5v_5\}$ with the condition (\ref{case211})
are equivalent to $\{v_5,v_3\}$ since the solution for Eqs.(\ref{Reqs}) is
\begin{equation}
k'_1=a_5e^{-\epsilon_4},\quad
k'_2=a_3+a_5\epsilon_1,\quad
k'_3=b_5e^{-\epsilon_4},\quad
k'_4=b_3+b_5\epsilon_1,\quad
\epsilon_2=\frac{a_1}{2a_5}.
\end{equation}

\subsubsection{the case of $\lambda\neq0 $ in the determined equations}

\textbf{\emph{Case 3:}} Not all $a_2, a_4$ and $a_6$ are zeroes. Without loss of generality, we also adopt $a_6\neq0$.

For $\lambda\neq0$ and $a_6\neq0$, by solving  Eqs.(\ref{fengbi}), we find
the relations
\begin{equation}
\label{relations}
a_1=\frac{a_4a_5}{2a_6},\quad a_2=\frac{a^2_4}{4a_6},\quad
a_3=-\frac{a^2_5+2a_4a_6}{4a_6}, \quad
b_1=\frac{a_4a_6b_5-a_5(a_4b_6-b_4a_6)}{2a^2_6},\quad
b_2=\frac{a_4(2a_6b_4-a_4b_6)}{4a^2_6}.
\end{equation}
Substituting (\ref{relations}) into Eqs.(\ref{8Eqs}), it leads to
an invariant for $\{w_1,w_2\}$:
\begin{equation}
\phi=\Delta_1 \equiv \frac{2(b_4+2b_3)a^2_6-a_5(a_5b_6-2b_5a_6)}{a_6(a_4b_6-a_6b_4)}.
\end{equation}
In condition of (\ref{relations}) and $\Delta_1=c$, choose the corresponding
representative element $\{v_6,(-\frac{1}{2}-\frac{c}{4})v_3+v_4\}$ and Eqs.(\ref{Reqs})have the solution
\begin{equation}
\begin{aligned}
&k'_1=a_6e^{-2\epsilon_4}, \quad k'_2=0, \quad
k'_3=2\epsilon_6b_4+b_6e^{-2\epsilon_4}-\frac{2\epsilon_6a_4b_6}{a_6},\quad
k'_4=\frac{a_6b_4-b_6a_4}{a_6},\\
&\epsilon_1=\frac{a_5}{2a_6},\quad
\epsilon_2=\frac{a_4}{4a_6},\quad
\epsilon_5=\frac{a_6b_5-a_5b_6}{a_4b_6-a_6b_4}e^{-\epsilon_4}.
\end{aligned}
\end{equation}

\textbf{\emph{Case 4:}} $a_2=a_4=a_6=0$.

For the same reason, we just think about
$a_5\neq0$. Taking $a_2=a_4=a_6=0$ with $a_5\neq0$ and $\lambda\neq0$ into (\ref{fengbi}), we have
\begin{equation}
b_1=\frac{a_1(2a_3b_6+a_5b_5)-a_3a_5b_4}{a^2_5},\quad
b_2=\frac{a_1(a_5b_4-a_1b_6)}{a^2_5},\quad
\lambda=\frac{2a_1b_6-a_5b_4}{a_5}\neq0.
\label{bbb}
\end{equation}
Substituting the relations (\ref{bbb}) into Eqs.(\ref{8Eqs}), we obtain
an invariant for $\{w_1,w_2\}$, i.e.
\begin{equation}
\phi=\Delta_2=\frac{4a_3(a_3b_6+a_5b_5)-2a^2_5(b_4+2b_3)}{a_5(2a_1b_6-a_5b_4)}.
\end{equation}
In this case, choose a representative element $\{v_5, (\frac{c}{4}-\frac{1}{2})v_3+v_4)\}$ for $\Delta_2=c$. Then solving
Eqs.(\ref{Reqs}), one get
\begin{equation}
\begin{aligned}
&k'_1=a_5e^{-\epsilon_4},\quad
k'_2=0,\quad
k'_3=\frac{(a_5b_5+2a_3b_6)e^{-\epsilon_4}+\epsilon_5(a_5b_4-2a_1b_6)}{a_5},\quad
k'_4=\frac{a_5b_4-2a_1b_6}{a_5},\\
&\epsilon_1=-\frac{a_3}{a_5},\quad
\epsilon_2=\frac{a_1}{2a_5},\quad
\epsilon_6=\frac{a_5b_6}{2(2a_1b_6-a_5b_4)}e^{-2\epsilon_4}.
\end{aligned}
\end{equation}

In summary, we have completed the construction of the two-dimensional
optimal system $\Theta_2$:
\begin{equation}
\label{optimal2}
\begin{aligned}
&\mathfrak{g}_1=\{v_6,v_3\}, \quad \mathfrak{g}_2=\{v_3+v_6,v_5\}, \quad \mathfrak{g}_2=\{-v_3+v_6,v_5\},\\
&\mathfrak{g}_4=\{v_6,v_5\}, \quad \mathfrak{g}_5=\{v_2+v_6,v_3\}, \quad \mathfrak{g}_6=\{-v_2+v_6,v_3\},\\
&\mathfrak{g}_7=\{v_1+v_6,v_3\},\quad \mathfrak{g}_8=\{-v_1+v_6,v_3\},\quad
\mathfrak{g}_{9}=\{v_5,v_3\}, \\
&\mathfrak{g}_{10}=\{v_6,v_4+\beta v_3\},\quad
\mathfrak{g}_{11}=\{v_5,v_4+\beta v_3\}.\quad
(\beta \in \mathbb{R})
\end{aligned}
\end{equation}

\textbf{\emph{Remark 4}}: The process of the construction ensures that all $\mathfrak{g}_i (i=1\cdots 11)$ are mutually inequivalent since each case
is closed. One can also easily find this inequivalence from the incompatibility of Eqs.(\ref{Reqs}).
In \cite{Qu1}, Chou et ale. gave a two-parameter optimal system $\{M_i\}^{10}_{1}$
for the same Lie algebra (\ref{vector}) of the heat equation and  showed their
inequivalences by sufficient  numerical invariants. One can see that $\{M_i\}^{10}_{1}$
in \cite{Qu1} are equivalent to our $\{\mathfrak{g}_{11},\mathfrak{g}_{10},\mathfrak{g}_{6},
\mathfrak{g}_{4},\mathfrak{g}_{1},\mathfrak{g}_{2},\mathfrak{g}_{3},
\mathfrak{g}_{8}, \mathfrak{g}_{5}, \mathfrak{g}_{9}
\}$, respectively. Furthermore, we realize that $\mathfrak{g}_{7}$ is
inequivalent to any of the elements in $\{M_i\}^{10}_{1}$. Hence, here the two-dimensional optimal system $\Theta_2$ given by (\ref{optimal2}) is
complete and  really optimal.

\section{the new approach for  the the (2+1)-dimensional
Navier-Stokes equation}

One of the most important open problems in fluid is the existence
and smoothness problem of the Navier-Stokes (NS) equation, which has been
recognized as the basic equation and the very starting point of all
problems in fluid physics~\cite{NS1,NS2}. In Ref.~\cite{NShu},
by means of the classical Lie symmetry method, we investigate
the (2+1)-dimensional Navier-Stokes equations:
\begin{equation}
\label{eq58}
\begin{aligned}
&\omega=\psi_{xx}+\psi_{yy},\\
&\omega_{t}+\psi_{x}\omega_{y}-\psi_{y}\omega_{x}-\gamma(\omega_{xx}+\omega_{yy})=0.
\end{aligned}
\end{equation}
One can rewrite Eqs.(\ref{eq58}) into
\begin{equation}
\label{NS}
\psi_{xxt}+\psi_{yyt}+\psi_{x}\psi_{xxy}+\psi_{x}\psi_{yyy}-\psi_{y}\psi_{xxx}-\psi_{y}\psi_{xyy}
 -\gamma(\psi_{xxxx}+2\psi_{xxyy}+\psi_{yyyy})=0.
\end{equation}
The associated vector fields for the one-parameter Lie group of the NS equation (\ref{NS}) are  given by
\begin{equation}
\begin{aligned}
\label{NSalgebra}
&v_1=\frac{x}{2}\partial_x+\frac{y}{2}\partial_y+t\partial_t, \quad
v_2=\partial_t, \quad
v_3=-yt\partial_x+xt\partial_y+\frac{x^2+y^2}{2}\partial_\psi, \quad
v_4=-y\partial_x+x\partial_y,\\
&v_5=f(t)\partial_x-f'(t)y\partial_\psi, \quad
 v_6=g(t)\partial_y+g'(t)x\partial_\psi, \quad
 v_7=h(t)\partial_\psi.
\end{aligned}
\end{equation}

Here, ignoring the discussion of the infinite dimensional
subalgebra, we apply the new approach to construct the two-dimensional
optimal system  and  the corresponding invariant solutions for the four-dimensional Lie algebra spanned by $v_1, v_2, v_3, v_4$ in (\ref{NSalgebra}).

The  commutator table and the adjoint representation
table for $\{v_1, v_2, v_3, v_4\}$ are given in table 3 and table 4,
respectively.
\begin{table}[htbp]
\centering
\caption{\label{table3}the commutator table for $\{v_1, v_2, v_3, v_4\}$}
\begin{tabular}{c|cccc}
\hline
     & $v_1$ &  $v_2$ & $ v_3$ &  $v_4$  \\
\hline
 $v_1$ & 0      &  $-v_2$ & $v_3$ &  $0$  \\
 \hline
 $v_2$ & $v_2$  & 0 &  $v_4$ & $0$  \\
 \hline
 $v_3$ & $-v_3$ & $-v_4$ & 0 & 0  \\
 \hline
 $v_4$ & 0 &  0 & 0 & 0  \\
 \hline
\end{tabular}
\end{table}

\begin{table}[htbp]
\centering
\caption{\label{table4}the adjoint representation table for $\{v_1, v_2, v_3, v_4\}$}
\begin{tabular}{c|cccc}
\hline
   Ad  & $v_1$ &  $v_2$ & $ v_3$ &  $v_4$  \\
\hline
 $v_1$ & $v_1$ & $e^{\epsilon} v_2$ & $e^{-\epsilon}v_3$ &  $v_4$  \\
 \hline
 $v_2$ & $v_1-\epsilon v_2$  & $v_2$ &  $v_3-\epsilon v_4$ & $v_4$  \\
 \hline
 $v_3$ & $v_1+\epsilon v_3$ & $v_2+\epsilon v_4$ & $v_3$ & $v_4 $ \\
 \hline
 $v_4$ & $v_1$ &  $v_2$ & $v_3$ & $v_4$  \\
\hline
\end{tabular}
\end{table}

\subsection{the  adjoint  transformation  matrix and the invariants equation}

Applying the adjoint action of $v_1$  to $w_1=
\sum\limits_{i=1} ^{4} a_iv_i$, we have
\begin{equation}
Ad_{exp(\epsilon_1 v_1)}(a_1v_1+a_2v_2+a_3v_3+a_4v_4)
=a_1v_1+a_2 e^{\epsilon_1}v_2+a_3e^{-\epsilon_1}v_3+a_4v_4.
\end{equation}
Hence the corresponding adjoint  transformation  matrix $A_1$ is saying
\begin{equation}
A_1=\left(
\begin{array}{cccc}
 1 & 0 & 0 & 0 \\
 0 & e^{\epsilon_1} & 0 & 0\\
 0 & 0 & e^{-\epsilon_1} & 0  \\
 0 & 0 & 0 & 1
\end{array}
\right).
\end{equation}
Similarly, one can get

 \begin{equation}
A_2=\left(
\begin{array}{cccc}
 1 & -\epsilon_2 & 0 & 0  \\
 0 & 1 & 0 & 0 \\
 0 & 0 & 1 & -\epsilon_2  \\
 0 &0 & 0 & 1
\end{array}
\right),\quad
A_3=\left(
\begin{array}{cccc}
 1 & 0 & \epsilon_3 & 0  \\
 0 & 1 & 0 & \epsilon_3 \\
 0 & 0 & 1 & 0 \\
 0 & 0 & 0 & 1
\end{array}
\right),
\quad
A_4=\left(
\begin{array}{cccc}
 1 & 0 & 0 & 0  \\
 0 & 1 & 0 & 0 \\
 0 & 0 & 1 & 0 \\
 0 & 0 & 0 & 1
\end{array}
\right)
\end{equation}
Then the most general adjoint  matrix $A$ can be taken as
\begin{equation}
A=A_1A_2A_3A_4=
\left(
\begin{array}{cccc}
 1 & -\epsilon_2 & \epsilon_3 & -\epsilon_2\epsilon_3 \\
 0 & e^{\epsilon_1} & 0 & e^{\epsilon_1}\epsilon_3\\
 0 & 0 & e^{-\epsilon_1} & -e^{-\epsilon_1}\epsilon_2 \\
 0 & 0 & 0 & 1
\end{array}
\right).
\end{equation}
Let
\begin{equation}
\label{NS12}
w_1=\sum\limits_{i=1} ^{4} a_iv_i, \quad w_2=\sum\limits_{j=1} ^{4} b_jv_j,
\end{equation}
For the general two-dimensional subalgebra $\{w_1, w_2\}$, the corresponding
invariant $\phi$ is some eight-dimensional function of $(a_1,\cdots, a_4, b_1, \cdots, b_4)$. Let $v=\sum\limits_{k=1} ^{4} c_kv_k$ be a general element of $\mathcal{G}$, then in conjunction with the commutator table 3, we have
\begin{equation}
\begin{aligned}
Ad_{g}(w_1)&=Ad_{exp(\epsilon v)}(w_1)\\
&=w_1-\epsilon [v,w_1]+ \frac{1}{2!}\epsilon^2[v,[v,w_1]]-\cdots \label{NSad1}
\\
&=(a_1v_1+\cdots+a_4v_4)-\epsilon [c_1v_1+\cdots+c_4v_4,a_1v_1+\cdots+a_4v_4]+o(\epsilon) \\
&=a_1v_1+(a_2-\epsilon(c_2a_1-c_1a_2))v_2+(a_3-\epsilon(c_1a_3-c_3a_1))v_3
+(a_4-\epsilon(c_2a_3-c_3a_2))v_4+o(\epsilon).
\end{aligned}
\end{equation}
Similarly, applying  the same adjoint action $v=\sum\limits_{k=1}^{4} c_kv_k$ to $w_2$, we get
\begin{equation}
\begin{aligned}
\label{NSad2}
Ad_{g}(w_2)=b_1v_1+(a_2-\epsilon(c_2b_1-c_1b_2))v_2+(b_3-\epsilon(c_1b_3-c_3b_1))v_3
+(b_4-\epsilon(c_2b_3-c_3b_2))v_4+o(\epsilon).
\end{aligned}
\end{equation}
Substituting (\ref{NSad1}) and (\ref{NSad2}) into Eq.(\ref{invariant}), we  treat two cases in the follows.

(a) When $\lambda=0$, the invariant function $\phi=\phi(a_1,a_2,a_3,a_4,b_1,b_2,b_3,b_4)$ is determined by

\begin{equation}
\label{NS8Eqs}
\begin{aligned}
a_1\phi_{a_1}+a_2\phi_{a_2}+a_3\phi_{a_3}+a_4\phi_{a_4}=0,\quad
a_1\phi_{b_1}+a_2\phi_{b_2}+a_3\phi_{b_3}+a_4\phi_{b_4}=0,\\
b_1\phi_{b_1}+b_2\phi_{b_2}+b_3\phi_{b_3}+b_4\phi_{b_4}=0,\quad
a_2\phi_{a_2}-a_3\phi_{a_3}+b_2\phi_{b_2}-b_3\phi_{b_3}=0,\\
a_1\phi_{a_2}+a_3\phi_{a_4}+b_1\phi_{b_2}+b_3\phi_{b_4}=0,\quad
a_1\phi_{a_3}+a_2\phi_{a_4}+b_1\phi_{b_3}+b_2\phi_{b_4}=0,\\
\end{aligned}
\end{equation}
and
\begin{equation}
b_1\phi_{a_1}+b_2\phi_{a_2}+b_3\phi_{a_3}+b_4\phi_{a_4}=0.
\label{NS1Eq}
\end{equation}

(b) When $\lambda \neq 0$, the invariant function $\phi=\phi(a_1,a_2,a_3,a_4,b_1,b_2,b_3,b_4)$ only  needs to meet
Eqs.(\ref{NS8Eqs}) in terms of  the condition $a_{12}=0$.

Now we set out to construct the two-dimensional optimal system for $\{v_1,v_2,v_3,v_4\}$.

\subsection{the construction of  two-dimensional optimal system for the heat equation}

Substituting (\ref{NS12}) into $[w_1,w_2]=\lambda w_1$,
the determined equations
are found as follows
\begin{equation}
\label{NSfengbi}
\begin{aligned}
 \lambda a_1=0, \quad
a_2b_1-a_1b_2=\lambda a_2, \quad
a_1b_3-a_3b_1=\lambda a_3, \quad
a_2b_3-a_3b_2=\lambda a_4.
\end{aligned}
\end{equation}

\subsubsection{the case of $\lambda=0$ in the determined equations (\ref{NSfengbi})}

\textbf{\emph{Case 1:}} Not all $a_1$ and $b_1$ are zeroes. Without loss of generality, we adopt $a_1\neq0$. Solving (\ref{NSfengbi}), we get

\begin{equation}
\label{NSeqi}
b_2=\frac{a_2b_1}{a_1},\quad b_3=\frac{a_3b_1}{a_1},
\end{equation}
with $a_1,a_2,a_3,a_4,b_1$ and  $b_4$ being arbitrary.

Substituting the  condition (\ref{NSeqi}) into Eqs.(\ref{NS8Eqs}) and Eq.(\ref{NS1Eq}), we find that $\phi= constant $.  Hence according to
(\ref{NSeqi}), select the corresponding representative element
$\{v_1,v_4\}$. Sinece  Eqs.(\ref{Reqs}) have the solution
\begin{equation}
k'_1=a_1, \quad k'_2=\frac{a_1a_4-a_2a_3}{a_1},\quad
k'_3=b_1, \quad k'_4=\frac{b_4a^2_1-a_2a_3b_1}{a^2_1},\quad
\epsilon_2=\frac{e^{\epsilon_1}a_2}{a_1}, \quad
\epsilon_3=-\frac{a_3}{a_1e^{\epsilon_1}},
\end{equation}
case (1) is equivalent to $\{v_1,v_4\}$.

\textbf{\emph{Case 2:}} $a_1=b_1=0$. Now the determined equations
(\ref{NSfengbi}) become
\begin{equation}
\label{NScase21}
a_2b_3-a_3b_2=0.
\end{equation}

\textbf{\emph{Case 2.1:}} Not all $a_2$ and $b_2$ are zeroes and we let
$a_2\neq 0$. From Eq.(\ref{NScase21}), we get $b_3=\frac{a_3b_2}{a_2}$.
Then by solving Eqs.(\ref{NS8Eqs}) and Eq.(\ref{NS1Eq}), we find  $\phi\equiv constant $.
In this case, there exist three circumstances in terms of the sign of $a_2a_3$.\\
(i). When $a_3=0$,  choose the  representative element
$\{v_2,v_4\}$ and  Eqs.(\ref{Reqs}) have the solution
\begin{equation}
k'_1=e^{\epsilon_1}a_2, \quad k'_2=e^{\epsilon_1}\epsilon_3a_2+a_4,
\quad k'_3=e^{\epsilon_1}b_2,\quad k'_4=e^{\epsilon_1}\epsilon_3b_2+b_4.
\end{equation}
(ii). For $a_2a_3>0$, we select $\{v_2+v_3,v_4\}$  as a  representative element.
Eqs.(\ref{Reqs}) hold for
\begin{equation}
k'_1=\sqrt{\frac{a_3}{a_2}}a_2, \quad k'_2=\sqrt{\frac{a_3}{a_2}}(\epsilon_3-\epsilon_2)a_2+a_4,
\quad k'_3=\sqrt{\frac{a_3}{a_2}}b_2,\quad k'_4=\sqrt{\frac{a_3}{a_2}}(\epsilon_3-\epsilon_2)b_2+b_4.
\end{equation}
(iii). For $a_2a_3<0$, we select $\{v_2-v_3,v_4\}$  as a  representative element.
and Eqs.(\ref{Reqs}) have the solution
\begin{equation}
k'_1=\sqrt{-\frac{a_3}{a_2}}a_2, \quad k'_2=\sqrt{-\frac{a_3}{a_2}}(\epsilon_3+\epsilon_2)a_2+a_4,
\quad k'_3=\sqrt{-\frac{a_3}{a_2}}b_2,\quad k'_4=\sqrt{\frac{a_3}{a_2}}(\epsilon_3+\epsilon_2)b_2+b_4.
\end{equation}

\textbf{\emph{Case 2.2:}} For $a_2=b_2=0$, Eqs.(\ref{Reqs}) always stand up and
the general two-dimensional Lie algebra becomes $\{a_3v_3+a_4v_4,b_3v_3+b_4v_4\}$. Then
if not all $a_3$ and $b_3$ are zeroes (and let
$a_3\neq 0$), it will equivalent to $\{v_3,v_4\}$ since that Eqs.(\ref{Reqs})
have the solution
\begin{equation}
k'_1=e^{-\epsilon_1}a_3, \quad k'_2=-e^{-\epsilon_1}\epsilon_2a_3+a_4,
\quad k_3=e^{-\epsilon_1}b_3, \quad k'_4=-e^{-\epsilon_1}\epsilon_2b_3+b_4.
\end{equation}

For the case of $a_3=b_3=0$, the  general two-dimensional Lie algebra
$\{a_4v_4,b_4v_4\}$ is trivial.

\subsubsection{the case of $\lambda \neq 0$ in the determined equations (\ref{NSfengbi})}

Substituting $\lambda \neq 0$ into Eqs.(\ref{NSfengbi}), there must be
$a_1=0$.

\textbf{\emph{Case 3:}} $a_2\neq 0$. Now, Eqs.(\ref{NSfengbi}) require
\begin{equation}
\label{NSrelations}
\lambda=b_1 (\neq 0), \quad a_3=0, \quad a_4=\frac{a_2b_3}{b_1}.
\end{equation}
Substituting (\ref{NSrelations}) into Eqs.(\ref{NS8Eqs}), it leads to
an invariant for $\{w_1,w_2\}$:
\begin{equation}
\phi=\Delta_3 \equiv \frac{b_1b_4-b_2b_3}{b^2_1}.
\end{equation}
In condition of (\ref{NSrelations}) and $\Delta_3=c$, choose the corresponding
representative element $\{v_2,v_1+cv_4\}$ and Eqs.(\ref{Reqs}) have the solution
\begin{equation}
k'_1=a_2e^{\epsilon_1}, \quad k'_2=0, \quad
k'_3=-\epsilon_2b_1+e^{\epsilon_1}b_2,\quad
k'_4=b_1,\quad \epsilon_3=-\frac{b_3}{b_1}e^{-\epsilon_1}.
\end{equation}

\textbf{\emph{Case 4:}} $a_2=0$. By solving Eqs.(\ref{NSfengbi}), we get
\begin{equation}
\label{NNSrelations}
\lambda=-b_1 (\neq 0),  \quad b_2=\frac{b_1a_4}{a_3}.
\end{equation}
Substituting (\ref{NNSrelations}) with $a_1=a_2=0$ into Eqs.(\ref{NS8Eqs}), one can obtain an invariant as follows:
\begin{equation}
\phi=\Delta_4 \equiv \frac{b_1b_4-b_2b_3}{b^2_1}.
\end{equation}
In condition of (\ref{NNSrelations}) and $\Delta_4=c$, select a
representative element $\{v_3,v_1+cv_4\}$ and Eqs.(\ref{Reqs}) have the solution
\begin{equation}
k'_1=a_3e^{-\epsilon_1}, \quad k'_2=0, \quad
k'_3=\epsilon_3b_1+e^{-\epsilon_1}b_3,\quad
k'_4=b_1,\quad \epsilon_2=\frac{a_4}{a_3}e^{\epsilon_1}.
\end{equation}

In summary, a two-dimensional optimal system $\Theta_2$ for
the four-dimensional Lie algebra spanned by $v_1, v_2, v_3, v_4$ in (\ref{NSalgebra}) is shown as follows:
\begin{equation}
\label{NSoptimal2}
\begin{aligned}
&\mathfrak{g}'_1=\{v_1,v_4\}, \quad \mathfrak{g}'_2=\{v_2,v_4\}, \quad \mathfrak{g}'_3=\{v_2+v_3,v_4\},\quad
\mathfrak{g}'_4=\{v_2-v_3,v_4\}, \\
&\mathfrak{g}'_5=\{v_3,v_4\}, \quad \mathfrak{g}'_6=\{v_2,v_1+cv_4\},\quad
\mathfrak{g}'_7=\{v_3,v_1+cv_4\}.\quad
(c \in \mathbb{R})
\end{aligned}
\end{equation}

\subsection{two-dimensional reductions for the  NS equation (\ref{NS})}
Using a two-dimensional Lie algebra, one can reduce the (2+1)-dimensional NS equation
to the ordinary equation and further get the invariant solutions. For each
two-parameter element in the two-dimensional optimal system (\ref{NSoptimal2}), there will be a corresponding class of group
invariant solutions which will ne determined from a reduced
ordinary differential equation. For the case of $\mathfrak{g}'_1=\{v_1,v_4\}$ and $\mathfrak{g}'_2=\{v_2,v_4\}$, one can refer to \cite{NShu}. The case of
$\mathfrak{g}'_5$ leads to no group invariant solutions. Then we just consider the rest elements in (\ref{NSoptimal2}).

(a) $\mathfrak{g}'_3=\{v_2+v_3,v_4\}$ and $\mathfrak{g}'_4=\{v_2-v_3,v_4\}$. By solving
$(v_2\pm v_3)(\psi)=0$ and  $v_4(\psi)=0$, we have $\psi=F(x^2+y^2)\pm \frac{1}{2}t(x^2+y^2)$. Substituting
it into the NS equation (\ref{NS}), one can get
\begin{equation}
\label{red1}
8\gamma[\xi^2F^{(4)}(\xi)+4\xi F^{'''}(\xi)+2F^{''}(\xi)]\mp1=0
\end{equation}
with $\xi=x^2+y^2$. By solving Eqs.~(\ref{red1}), we find
$\mathfrak{g}'_3$ and $\mathfrak{g}'_4$ lead to the same
group invariant solution
\begin{equation}
\psi=c_1+c_2(x^2+y^2)+c_3\ln(x^2+y^2)+c_4(x^2+y^2)
[\ln(x^2+y^2)-1]+\frac{1}{32\gamma}(x^2+y^2)^2+\frac{1}{2}t(x^2+y^2).
\end{equation}

(b) $\mathfrak{g}'_6=\{v_2,v_1+cv_4\}$. From $v_2(\psi)=0$ and $(v_1+cv_4)(\psi)=0$, one can get $\psi=F(\arctan(\frac{y}{x})-c\ln(x^2+y^2))$. Substituting
it into (\ref{NS}) and integrating the reduced equation
once, we have
\begin{equation}
\label{red2}
\gamma[(4c^2+1)G''(\xi)+8cG'(\xi)+4G(\xi)]-G^2(\xi)=0, \quad (G(\xi)=F'(\xi))
\end{equation}
with $\xi=\arctan(\frac{y}{x})-c\ln(x^2+y^2)$.
Specially, when $c=0$ in Eq.~(\ref{red2}), there is  a solution
\begin{equation}
G(\xi)=-6\gamma \rm{sech}^2(\xi+c_0)+4\gamma.
\end{equation}
Then it leads to the solution of the NS equation
\begin{equation}
\psi=-6\gamma \tanh (\arctan(\frac{y}{x})+c_0)+4\gamma \arctan(\frac{y}{x})+c_1.
\end{equation}

(c) $\mathfrak{g}'_7=\{v_3,v_1+cv_4\}$. In this case, we have
$\psi=\arctan(\frac{x}{y})+c\ln(t)+F(\frac{x^2+y^2}{t})$. The reduced equation for Eqe~(\ref{NS}) is
\begin{equation}
4\gamma Z^2 G''(Z)+Z(Z+8\gamma-2)G'(Z)+(Z-2)G(Z)=0, \quad (F'(Z)=G(Z))
\end{equation}
with $Z=\frac{x^2+y^2}{t}$.
In particular, for $\gamma=\frac{1}{4}$, we obtain a solution
\begin{equation}
\psi=\arctan(\frac{x}{y})+c\ln(t)+c_1+c_2\ln(Z)+c_3(2\ln(Z)+3e^{-Z}+Ze^{-Z}
+2\rm{Ei}(1,Z)+\frac{4}{3});
\end{equation}
for $\gamma=1$, there is
\begin{equation}
\psi=\arctan(\frac{x}{y})+c\ln(t)+c_1+c_2\ln(Z)
+\frac{105}{8}c_3 \bigg(-\sqrt{\pi
}\rm{erf}(\frac{\sqrt{Z}}{2})+
\sqrt{Z}
\rm{hypergeom}([\frac{1}{2},\frac{1}{2}],[\frac{3}{2},\frac{3}{2}],-\frac{Z}{4})\bigg).
\end{equation}

\section{Summary and discussion}\label{conclusion}

Since many important equations arising from
physics are of low dimensions, only
the determination of small parameter optimal systems can
reduce them to ODEs which often lead to inequivalent
group invariant solutions. In this paper, we give an elementary algorithm for constructing two-dimensional optimal system which only
depends on fragments of the theory of Lie algebras. The  intrinsic
idea of our method is that every element in the optimal system
corresponds to different values of invariants, the  definition of which
have been refined in this paper.
Thanks to these  invariants which are often overlooked except the Killing form in the almost existing methods, all the elements  in the two-dimensional optimal system are found one by one and their  inequivalences are evident with no further proof. Moreover,
the construction of two-dimensional optimal system by us just starts
from subalgebras, and does not require the  prior one-dimensional optimal system as usual.

Before manipulating the given algorithm to construct two-dimensional optimal systems, one should compute two-dimensional subalgebras with ``the determined equations", the general adjoint transformation matrix and the invariants equations, which seem much complicated but can all be  carried out in mechanization with the compute software \emph{``Maple"}.
A new method is shown to provide all the invariants for the  two-dimensional subalgebras, which is based on
the idea of ``invariant" under the meaning of adjoint transformation and
combination act. Applying  the algorithm to (1+1)-dimensional heat equation and (2+1)-dimensional NS equation,  we obtain their
two-dimensional optimal systems respectively. For the heat equation,
the obtained two-dimensional optimal system contains eleven elements,
which are discovered  more comprehensive than that in Ref.~\cite{Qu1} after a detailed comparison. For the (2+1)-dimensional NS equation, all
the reduced ordinary equations and some exact group invariant solutions which come from the obtained two-dimensional optimal system are found.
Since group invariant solutions are connected with the specific differential systems, inequivalent subalgebras may lead to the same solutions.

The algorithm considered in this paper is elementary and practical, without too much algebraic knowledge. Due to the programmatic process,
to give a Maple package on the computer for two-dimensional optimal system is necessary and under our consideration. How to apply all the invariants to construct r-parameter ($r>2$) optimal systems is also an
interesting job.

\section*{Acknowledgments}

This work is supported by Zhejiang Provincial Natural Science Foundation of China under Grant No. LQ13A010014, the National Natural Science
Foundation of China (Grant Nos. 11326164, 10914294, 11075055 and 11275072), the Research Fund for the Doctoral Program of Higher Education of China (Grant No. 20120076110024).

\bibliography{<your-bib-database>}

\begin{thebibliography}{00}
\bibitem{Ovsiannikov} L. V. Ovsiannikov, Group Analysis of Differential Equations, Academic Press, New York, 1982.
\bibitem{Galas} Galas, F. and Richet, E. W.: Physica D 50 (1991), 297¨C307.
\bibitem{Ibragimov} N. H. Ibragimov, Lie Group Analysis of Differential Equations, Vol. 1, CRC Press, Boca Raton, FL, 1994.

\bibitem{Olver} P. J. Olver, Applications of Lie groups to differential equations, 2nd ed. Springer, New York, 1993.

\bibitem{Patera1} J. Patera, P. Winternitz, and H. Zassenhaus, Continuous subgroups of the fundamental groups of physics. I General method and the Poincar\'{e} group, J. Math. Phys. 16 (1975), 1597-1614.
\bibitem{Patera2} J. Patera, R. T. Sharp, P. Winternitz, and H. Zassenhaus, Invariants of real low dimension Lie algebras, J. Math. Phys. 17 (1976), 986-994.
\bibitem{Patera3} J. Patera, R. T. Sharp, P. Winternitz, and H. Zassenhaus, Subgroups of the Poincar\'{e} group and their invariants, J. Math. Phys. 17 (1976), 977¨C985.

\bibitem{Patera4} J. Patera and P. Winternitz, Subalgebras of real three- and four-dimensional Lie algebras, 18 (1977) 1449-1455.





\bibitem{Coggeshalla} S. V. Coggeshalla and J. Meyer-ter-Vehn, J. Math. Phys. 33 3585 (1992).

\bibitem{Qu1} K. S. Chou, G. X. Li, and C. Z. Qu, A Note on Optimal Systems for the Heat Equation, Journal of Mathematical Analysis and Applications, 261, 741-751 (2001).

\bibitem{Qu2} K. S. Chou and C. Z. Qu, Optimal Systems and Group Classification of (1+2)-Dimensional Heat Equation, Acta Applicandae Mathematicae, 83, 257-287, 2004.
\bibitem{Qu3} K. S. Chou and G. X. Li, Optimal Systems and Invariant solutions for the Curve Shorting Problem, communications in Analysis and Geometry, 10 (2), 241-274, 2002.



\bibitem{NS1} D. Sundkvist, V. Krasnoselskikh, P.K. Shukla,
A. Vaivads, M. Andr, S. Buchert, and H. Rme, Nature
(London) 436, 825 (2005);
\bibitem{NS2} G. Pedrizzetti, Phys. Rev.
Lett. 94, 194502 (2005).
\bibitem{NShu} X. R. Hu, Z. Z. Dong, F. Huang and Y. Chen, Z. Naturforsch. 65a, 504-510 (2010);


\end{thebibliography}

\end{document}